# Equilibrium condition in Insurance Pricing: a particular case

Renato Ghisellini, Sistemi Srl, Piazza 4 Novembre 1, Milano, Italia, e-mail: sistemi.srl@galactica.it

*Acknowledgements*: Special thanks to Professor Neil Doherty for very helpful discussions and suggestions.

**Abstract**

Equilibrium pricing has been proven to underlie the Insured rational expectancy of premia additivity for composition of policies fully covering independent risks.

**Introduction**

In a previus paper (Ghisellini (1998)), the form of the acceptable disutility function for a rational Insured has been singled out by considering:

a) the Insured rational expectancy of premia additivity in case policies fully covering independent risks are joined together,

b) the fact that, in the risk aversion hypothesis, the only disutility function family coherent with both the constraint contained in point a) and the (possible) policy fairness is given by

$$U(l) = \rho(e^{\frac{l}{\rho}} - 1), \qquad (1)$$

where $l>0$ represents a cost (for instance the total loss amount) and $\rho>0$ is the parameter known in literature as "risk tolerance".

The point a) results from the fact that the coverage configurations considered - the first one corresponding to the joint action of two separated policies $C_1$ and $C_2$ fully covering two independent risk $K_1$ and $K_2$, the second corresponding to a single policy $C$ fully covering the composition $K$ of the risks $K_1$ and $K_2$ (the global risk) - are identical in their effects, being identical their total pay off (indemnity) distribution (it is worth noting that the present discussion is based on the further assumption - not specified in the above-mentioned paper - that the two independent risks $K_1$, $K_2$ and the coverages $C_1$, $C_2$ and $C$ exist).

For this reason, in the particular case here considered, the expectancy of the rational Insured about policy premia is expressed by $P = P_1 + P_2$, where $P$ is the premium of the policy $C$, $P_1$ and $P_2$ are the premia for policies $C_1$ and $C_2$.

The equilibrium condition underlying this argument constitutes the object of the present note.

### Premia additivity as equilibrium condition

Given the independent risks $K_1$ and $K_2$, let $\tilde{X}_1$ and $\tilde{X}_2$ be the respective stochastic functions describing the total loss values. Let $K$ be the composition of the two risks $K_1$ and $K_2$, $\tilde{X} = \tilde{X}_1 + \tilde{X}_2$ being the stochastic function describing the total loss value related to the global risk $K$.

Let us consider the two coverage configurations specified in the Introduction. In the first one, risks $K_1$ and $K_2$ are fully covered (retention identically zero) by the separated policies $C_1$ and $C_2$ having $\tilde{I}_1 = \tilde{X}_1$ and $\tilde{I}_2 = \tilde{X}_2$ as total indemnities, $P_1$ and $P_2$ being the premia. In the second configuration, the global risk $K$ is fully covered by the single policy $C$, the total indemnity being $\tilde{I} = \tilde{X}$ and the premium being $P$.

*Proposition:*

Under the hypothesis above specified, and under the further hypothesis that credit risk or other aspects linked to Insurer rating are negligible, if the Insured and the Insurers act in the Insurance market as rational players the additive composition of the premia

$$P = P_1 + P_2 \qquad (2)$$

corresponds to the equilibrium condition. Moreover, this equilibrium pricing condition is unique.

*Proof*

*1. Total pay-off equivalence:*

Being the stochastic functions $\tilde{X}_1$ and $\tilde{X}_2$ mutually independent, coverages $C_1$ and $C_2$ result mutually independent as well because $\tilde{I}_1 = \tilde{X}_1$ and $\tilde{I}_2 = \tilde{X}_2$. The combined action $\tilde{I}_1 + \tilde{I}_2$ of the two separated policies $C_1$ and $C_2$ is so associated to a distribution function that is decomposable in the following way:

$$D_{\tilde{I}_1 + \tilde{I}_2} = D_{\tilde{I}_1} * D_{\tilde{I}_2} \qquad (3)$$

where "$D$" means distribution of the indicated variable and "$*$" means convolution.

But,

$$D_{\tilde{I}_1} * D_{\tilde{I}_2} = D_{\tilde{X}_1} * D_{\tilde{X}_2} = D_{\tilde{X}} = D_{\tilde{I}} \qquad (4)$$

and so:

$$D_{\tilde{I}_1 + \tilde{I}_2} = D_{\tilde{I}}. \qquad (5)$$

Equation (5) demonstrates the fact that the total pay-off of the policy $C$ is completely equivalent to the one associated to combined action of the separated policies $C_1$ e $C_2$. For that, policy $C$ is completely equivalent to the two policies $C_1$ and $C_2$ taken together.

*2. Equilibrium*

Let

$$\hat{A} = \{A^1, A^2, ..., A^n\} \qquad (6)$$

be the set of all the Insurers acting on the market.

Let

$$\hat{C}_1 = \{C_1^1, C_1^2, ..., C_1^{v_1}\} \neq \phi, \quad (7)$$

$$\hat{C}_2 = \{C_2^1, C_2^2, ..., C_2^{v_2}\} \neq \phi \quad (8)$$

and

$$\hat{C} = \{C^1, C^2, ..., C^v\} \neq \phi \quad (9)$$

be respectively the sets of the full coverages offered by the market $\hat{A}$ for the risks $K_1$, $K_2$ and $K$ with $v_1 \leq n$, $v_2 \leq n$, $v \leq n$ and $\phi$ = the empty set. As specified in the Introduction and in (7) ÷ (9), it is assumed that for each considered risk the market $\hat{A}$ is able to offer at least one coverage.

Let also

$$\hat{P}_1 = \{P_1^1, P_1^2, ..., P_1^{v_1}\}, \quad (10)$$

$$\hat{P}_2 = \{P_2^1, P_2^2, ..., P_2^{v_2}\} \quad (11)$$

and

$$\hat{P} = \{P^1, P^2, ..., P^v\} \quad (12)$$

be the the sets of premia corresponding to $\hat{C}_1$, $\hat{C}_2$ and $\hat{C}$.

The hypothesis of rationality and the hypothesis of negligible credit risk together imply that the premia $P_1$, $P_2$ and $P$ of the policies $C_1$, $C_2$ and $C$ which the Insured would consider, necessarily fulfil the conditions:

$$P_1 = \min\{P_1^1, P_1^2, ..., P_1^{v_1}\}, \quad (13)$$

$$P_2 = \min\{P_2^1, P_2^2, ..., P_2^{v_2}\} \quad (14)$$

and

$$P = \min\{P^1, P^2, ..., P^v\}. \quad (15)$$

Let $A_1^{m_1}$, $A_2^{m_2}$ and $A^m$ be the Insurers supplying the coverages $C_1$, $C_2$ and $C$.

The fact that (2) expresses an equilibrium condition is easily seen by analyzing what happens in the cases $P<P_1+P_2$ and $P>P_1+P_2$.

If

$$P < P_1 + P_2, \tag{16}$$

since the equivalence condition (5), the Insured would prefer $C$ instead of the policies $C_1$ and $C_2$ taken together. Again since (5), the Insurer $A^m$ - if for instance requested by the Insured - could always be able to split the policy $C$ into two separated policies $\Gamma_1$ and $\Gamma_2$, fully covering respectively risks $K_1$ and $K_2$, and to supply them to the Insured at the premia $\pi_1$ and $\pi_2$ valued

$$\pi_1 + \pi_2 = P \leq P_1 + P_2 \tag{17}$$

(Borsh (1962)).

For this reason, at least one of the two following equation should hold true:

$$\pi_1 \leq P_1, \tag{18}$$

$$\pi_2 \leq P_2 \tag{19}$$

this fact contradicting either (13) or (14).

If

$$P > P_1 + P_2, \tag{20}$$

the Insurers $A_1^{m_1}$ and $A_2^{m_2}$ would act, exploiting mutual independence of the risks $K_1$ and $K_2$ and so eq. (5), in order to supply a policy $\Gamma$ equivalent to $C$ but cheaper than $C$, which would yield for both of them an economic result higher than the premia of policies $C_1$ and $C_2$ by themselves supplied. One of the possible strategies could be the one where $A_1^{m_1}$ supplies a full insurance of the global risk K for the premium

$$\Pi = \frac{P + P_1 + P_2}{2} < P, \tag{21}$$

the risk $K_2$ being fully reinsured by $A_2^{m_2}$ (which would take also client-service and administrative costs for risk $K_2$ as in the case of policy $C_2$) for the premium

$$\Pi_2 = P_2 + \frac{P - P_1 - P_2}{4} > P_2. \tag{22}$$

The risk assumed by $A_1^{m_1}$ would result $K_1$, the net premium being given by

$$\Pi_1 = \Pi - \Pi_2 = P_1 + \frac{P - P_1 - P_2}{4} > P_1. \tag{23}$$

Equation (21) would contradict equation (15).

These last facts, together with the fact that the sets (7)-(9) are not empty, allow saying that (2) constitute the equilibrium condition and that this equilibrium condition is unique.

*Q.E.D.*

It is worth noting that eq. (13) - (15) constitute a necessary but not sufficient conditions for the actual choice of the Insured. This means that the fact that the Insured is "free" to purchase or not the policies (depending on their values) is preserved and does not affect the present proof.

## Conclusions

The premia additivity for composition of policies fully covering independent risks does correspond to an equilibrium pricing condition. This objective condition, intrinsic in a rational Insurance market, defines the rational pricing expectancy of the Insured in the considered particular case and acts by constraining the Insured disutility function form as described in Ghisellini (1998).